\documentclass[a4paper,12pt]{article}

\usepackage[latin1]{inputenc}
\usepackage[english]{babel}
\usepackage{amssymb,amsfonts,mathrsfs}
\usepackage{amsmath}
\usepackage{amsthm}
\usepackage{amscd}
\usepackage{mathrsfs}
\usepackage{eufrak}
\usepackage{latexsym}
\usepackage{color}
\usepackage{enumerate}

\newtheorem{theo}{Theorem}[section]
\newtheorem{prop}[theo]{Proposition}
\newtheorem{lemma}[theo]{Lemma}
\newtheorem{coro}[theo]{Corollary}

\theoremstyle{definition}
\newtheorem{defi}[theo]{Definition}
\newtheorem{exa}[theo]{Example}
\newtheorem{oss}[theo]{Remark}

\author{Fabrizio Brienza\footnote{DISIM - Dipartimento di Ingegneria e Scienze 
dell'Informazione e Matematica
-- Universit\`a degli Studi dell'Aquila -- Via Vetoio snc, 67010
COPPITO (AQ) Italy -- e-mail:
\texttt{brienzafabrizio@yahoo.it}.}, 
Anna Guerrieri\footnote{DISIM - Dipartimento di Ingegneria e Scienze 
dell'Informazione e Matematica
-- Universit\`{a} degli Studi dell'Aquila -- Via Vetoio snc, 67010
COPPITO (AQ) Italy -- e-mail:
\texttt{guerran@univaq.it}.} 
}
\date{}

\title{\textbf{INITIAL IDEALS OF BOREL TYPE}}

\begin{document}

\maketitle

\begin{abstract} 
In this paper we use some results related to regularity,
Betti numbers and reduction of generic initial ideals,
showing their stability in passing from an ideal to 
its initial ideal if the last has some simple properties.
\end{abstract}

\section{Introduction}
Throughout the paper $R = K[x_1,\ldots,x_n]$ is
the polynomial ring in $n$ variables over an infinite field $K$,
$<$ a monomial order on R with $x_1 > x_2 > \cdots > x_n$
and $M$ a graded $R-$module. It is well know that for a graded 
ideal $I \subseteq R$ (an ideal generated by homogeneous elements)  
there exists a nonempty open set $U$ of linear automorphisms of R such
that in$_{<}(\alpha I)$ does not depend on $\alpha \in U$. The resulting 
initial ideal, gin$_<(I)$ is called the generic initial ideal of $I$ with respect 
to $<$. Generic initial ideals are Borel-fixed [\cite{bs}, Theorem 2.8] 
and are even strongly stable if the base field is of characteristic 0 [see \cite{bs}]. 
Passing to the generic initial ideal with reverse lexicographic order preserves the
extremal Betti numbers [\cite{bcp}, Theorem 1.6] and the reduction number [\cite{t}, Theorem $4.3$].  
However it is difficult to compute gin, because one does not have much information about the open
subset $U$, besides the fact that it is dense in $K^m$ in the standard topology ($m = n^2$) and 
therefore hard to avoid. 
Thus if we pick $x \in K^m$ randomly, i.e. ``generically enough'', 
then most likely $x$ will belong to $U$ and this is how most computer algebra 
systems compute gin$_<(I)$. An uncertainty though remains. In the wake of the works of
Bermejo and Gimenez \cite{bg}, Conca, Herzog and Hibi \cite{chh}, and Trung \cite{t}, we avoid gins to show the same results on numerical 
invariants for Borel type ideals. Bayer and Stillman in \cite{bs} prove that 
ideals Borel-fixed are of Borel type [for the definition see Section $2$] even
though the converse is clearly not true. 
A monomial ideal is of Borel type if and only if all the
annihilator modules associated to the sequence $\{x_n , x_{n-1} , \ldots , x_{1}\}$
are zero dimensional [Proposition \ref{ultima}]. Using this fact
we prove that the extremal Betti numbers [Theorem \ref{principale}] and
the reduction number with respect the sequence $\{x_n , \ldots , x_{n-d+1}\}$, where
$d =$ dim$(R/I)$ [see Section 4], 
are preserved in passing to the initial ideal if the latter is of Borel type. 
We also observe that the annihilator numbers are preserved 
[Corollary \ref{4}]. The annihilator numbers of a filter regular sequence are intimately related 
to the extremal Betti numbers, in the sense that the two diagrams are specular; in this way one attains also the 
information about extremal Betti numbers. We choose this approach because
the annihilator numbers are easy to compute, since they are in fact colons.
The present work is divided in four sections. In the first we recall some basic properties 
related to Borel-type ideals and to the annihilator numbers of a filter regular sequence. In the second section
we prove that the extremal Betti numbers and the annihilator
numbers of $I$ and in$_<(I)$ are equal in the case in$_<(I)$ is a Borel-type
ideal. Then we study the rigidity of resolutions of $I$ and in$_<(I)$, if in$_<(I)$ is of 
Borel-type and we show that, if $I$ is an ideal with an initial ideal of Borel-type, 
we don't have necessarily the rigidity of the resolution. 
In the third section we see that also the reduction numbers of $I$ and in$_<(I)$ with respect
the sequence $\{x_n, \ldots , x_{n-d+1}\}$
are the same, if the last ideal is of Borel-type. In the last section we compare the Borel-type ideals
to the quasi stable ideals of Hashemi, Schweinfurter and Seiler.
Recently Hashemi, Schweinfurter and Seiler in \cite{hss} proved, using the Pommaret bases 
(a special class of Gr\"obner basis with additional combinatorial properties), that, if $I$ has
a finite Pommaret basis, $I$ and in$(I)$ have the same extremal Betti numbers. 
When $I$ has a finite Pommaret basis, it is called a quasi-stable ideal. 
It is possible to prove that monomial quasi-stable ideals are mirror images of 
Borel type ideals. It is however possible to give examples of monomial quasi-stable 
ideals that are not Borel type and viceversa, even though a suitable change of variables transforms the one in
the others. Hashemi, Schweinfurter and Seiler also prove that given $I$
a polynomial ideal, there always exists a change of variables such that $I$
has a finite Pommaret basis. These variables are called $\delta-$regular. From their
result and the observation above, it follows that after a suitable change of variables, 
in$(I)$ becomes an ideal of Borel type.
It is now clear that, if $I$ has an initial ideal of Borel type, after a suitable
change of basis, it means that $I$ has a finite Pommaret basis.

\section{Preliminary notions and Borel type ideals}
Given the assumptions described in the introduction, we recall
some basic notions. Let
\[  
\cdots \longrightarrow F_i = \bigoplus_{j \geqslant 0}{R(- j)^{\beta_{ij}}}  
\longrightarrow \cdots \longrightarrow F_0 = \bigoplus_{j \geqslant 0}{R(-j)^{\beta_{0j}}} 
\longrightarrow M \longrightarrow 0,
\]
be the minimimal free resolution of $M$, where $\beta_{ij}(M) \geqslant 0$ is the rank of
the shift $ - j $ in $i-$th position. The minimum length of such a free resolution is called the
\textbf{projective dimension} of $M$ over $R$ and it is written pd$(M)$. $\beta_{ij}(M)$, 
for short $\beta_{ij}$, are called the \textbf{Betti graded numbers} of $M$. Betti numbers 
have been widely investigated and for the general theory we refer to \cite{bh2}. 
It is well known, [see \cite{bh}], that $\displaystyle{\beta_{ij} = \text{dim}_{k}{\text{Tor}_{i}^{R}(M,k)}_j = 
\text{dim}_{k}{H_{i}(\mathbb{F} \otimes k)}_j}$, and $\displaystyle{\text{pd}(M)} = \text{max} 
\{ i : \beta_{ij} \neq 0 \ \textrm{for some} \ j  \}.$

\begin{defi} \upshape Let $a_j$ be the maximum degree of 
the generators of $F_j$. Then $\text{reg}(M) = \max\{a_j - j : j>0\}$
is called the \textbf{Castenuovo-Mumford regularity}, or simply \textbf{regularity}, of $M$. 
\end{defi}
The regularity is an important invariant which measures the complexity of the 
given module; for the theory of regularity see \cite{e}. It is well known the 
connection with Betti numbers, in fact
\[
\text{reg}(M) = \max\{ j : \beta_{i,i+j}(M) \ne 0 \ \ \text{for some} \ i \}.
\] 
Let now $\mathbf{m}$ be the maximal graded ideal of $R$. Suppose that
$M$ is finitely generated, then we denote by $H_{\mathbf{m}}^i(M)$ the 
$i-$\textbf{th} \textbf{local cohomology} of $M$  with respect to 
$\mathbf{m}$ [see \cite{hh}]. Since $H_{\mathbf{m}}^i(M)$ is artinian, we may
consider the numbers $a_i(M) = \text{max}\{ n : (H_{\mathbf{m}}^i(M))_n \neq 0 \}$
assuming that $a_i(M) = - \infty$, if $H_{\mathbf{m}}^i(M) = 0$. 
The Castelnuovo-Mumford regularity of $M$ can be defined also as
\[
\text{reg}(M) = \text{max}\{ a_i(M) + i : i \geqslant 0 \}.
\]
Moreover the \textbf{largest non-vanishing degree} of local
cohomology modules is defined as the number
$a^{\ast}(M) = \text{max}\{ a_i(M) : i \geqslant 0 \}.$

The Castelnuovo-Mumford regularity and the largest non-vanishing
degree for local cohomology modules can be viewed as special cases
of the more general invariants:
\[
\text{reg}_t(M) = \text{max}\{ a_i(M) + i : i \leqslant t \}, 
\]
\[
a^{\ast}_t(M) = \text{max}\{ a_i(M) : i \leqslant t \},
\]
where $t \in \{0,\ldots,d\}$, where $d =$ dim$(M)$. These invariants have been studied in
\cite{t1}, \cite{t2}, \cite{t3}.

Now let $I \subseteq R$ be a graded ideal. We
define the ideals
\[
I : (x_1,\ldots,x_i)^\infty = 
\displaystyle\bigcup_{k \geqslant 0} I : (x_1,\ldots,x_i)^k, 
\] 
for $i=1,\ldots,n$. By abuse of notation $I : x_j^\infty = I : (x_j)^\infty$.

Accordingly with \cite{hh}, we give the following definition.
 
\begin{defi} \upshape
The ideal $I$ is said to be of \textbf{Borel type}, or a 
\textbf{weakly stable} ideal, if $ I : x_j^\infty = 
I : (x_1,\ldots,x_j)^\infty, \ \ \text{for all} \ j=1,\ldots,n.$
\end{defi}

We recall that an ideal $I \subseteq R$ is said to be Borel-fixed if
$\alpha(I) = I$ for all $\alpha \in \mathcal{B}$, where $\mathcal{B}$
is the Borel subgroup of GL$_n(K)$, that is the subgroup of all 
non-singular upper triangular matrices.

\begin{oss} \upshape Let $I$ be an ideal of $R$. An example
of ideal of Borel type is the generic initial ideal of $I$, gin$(I)$. 
This depends on the fact that gin$(I)$ is Borel-fixed as showed in
\cite{bs}.
\end{oss}

\begin{prop} \label{Bti} \emph{[\cite{hh}, Proposition 4.2.9]} 
Let $I \subseteq R$ be a graded monomial ideal,
that is an ideal generated by monomials. 
The following conditions are equivalent:
\begin{enumerate}
\item[1.]{$I$ is of Borel type;}

\item[2.]{for each monomial $u \in I$ and all integers $i,j,s$ with 
$1 \leqslant j < i \leqslant n$ and $s>0$ such that $x_i^s | u$, 
there exists an integer $t \geqslant 0$ such that $x_j^t(u/x_i^s) \in I$;}

\item[3.]{for each monomial $u \in I$ and all integers $i,j$ with 
$1 \leqslant j < i \leqslant n$, there exists an integer $t \geqslant 0$
such that $x_j^t(u/x_i^{\nu_i(u)}) \in I$, where $\nu_i(u)$ is the
highest power of $x_i$ which divides $u$;}

\item[4.]{if $p \in$ Ass(R/I), then $p = (x_1,\ldots,x_j)$ for some $j$.}
\end{enumerate}
\end{prop}

We show a class of ideals whose initial ideals
are of Borel type in any characteristic [see \cite{co0}]:
\begin{exa} \upshape Let $I = ((ax+by)^2 , (cx+dy)^2)$ be ideals
in $K[x,y]$ with $a,b,c,d \in K^{\ast}$. Let $x > y$ a monomial order. 
Then, using CoCoA, one can see that
\[ \text{in}_<(I) =
\begin{cases}
(x^2,xy,y^3) & \text{if char}(K) \neq 2 \ \text{and} \ ab \neq 0; \\
(x^2,y^2)  & \text{if char}(K) = 2 .
\end{cases}
\]
We can easily see that in$_<(I)$ is Borel type in both cases.
\end{exa}

\begin{defi} \upshape Let $I \subseteq R$ be a monomial ideal. Then $I$
is \textbf{strongly stable} if one has $x_i(u/x_j) \in I$ for all monomials $u \in I$
and all $i < j$ such that $x_j$ divides $u$.

\end{defi}

\begin{oss} \upshape It easily follows that a strongly stable monomial ideal is always
of Borel type, moreover a strongly stable ideal is Borel-fixed. Furthermore in a characteristic zero 
field an ideal is Borel-fixed if and only if it is strongly stable. 
[\cite{hh}, Proposition $4.2.4$].
\end{oss}

We use this remark to give an example of an ideal that is of Borel type
but not Borel-fixed.
\begin{exa} \upshape Let $I = (x_1^3 , x_1x_2^2, x_1^3x_2, x_1x_3^2)$ be an ideal
in $K[x_1,x_2,x_3]$ with char$(K) = 0$. Since $x_3 | x_1x_3^2$ but $x_1x_2x_3 \notin I$,
we know that $I$ is not Borel-fixed, since it is not a strongly stable ideal.
Checking the condition $(3)$ of Proposition \ref{Bti}, we show
that $I$ is of Borel type. Since $x_2 | x_1x_2^2$ and $\nu_2(x_1x_2^2) = 2$ we 
have to show that there exists an integer $t \geqslant 0$ such that
$x_1^{t+1} \in I$. It is sufficient to pick $t = 2$. For the next generator we have
that $x_2 | x_1^3x_2$ with $\nu_2(x_1^3x_2) = 1$ and so for $t = 0$ we know that
$x_1^{t+3} \in I$. For the last generator $x_1 | x_1x_3^2$ and $x_3 | x_1x_3^2$ with
$\nu_1(x_1x_3^2) = 1$ and $\nu_{3}(x_1x_3^2) = 2$. If we pick $t \geq 2$ we have that
$x_1x_2^t \in I$ and $x_1^{t+1} \in I$. Hence $I$ is of Borel type. 
\end{exa}

\begin{defi} \upshape Let $M$ be an $R$-module. 
An element $y \in R_1$ is said to be \textbf{filter regular} on $M$ if the 
multiplication map $y : M_{i-1} \to M_i$ is injective for all $i \gg 0$. 
The elements $y_1 , \ldots , y_r \in R_1$ form a \textbf{filter regular sequence} on $M$,
if $y_i$ is filter regular on $M/(y_1,\ldots,y_{i-1})M$, for all $i=1,\ldots,r$.
\end{defi}

\begin{oss} \upshape Immediately follows that, $y \in R_1$ is filter regular
on $M$, if and only if the ideal $(0:_{M} y)$ has finite length.
Herz\"og and Hibi in \cite{hh} show that if $|K| = \infty$ always exists a $K$-basis of $R_1$
that is a filter regular sequence on $M$. 
\end{oss}

\begin{defi} \upshape Let $M$ be a finitely generated graded $R$-module.
Let $\mathbf{y}=y_1,\ldots,y_n$ elements in $R_1$.
We denote by $A_{i-1}(\mathbf{y};M)$ the graded $R$-module
$\left( 0 :_{M/(y_1,\ldots,y_{i-1})M} y_i \right)$. The numbers
$\alpha_{ij}(\mathbf{y};M) = \dim_K A_i(\mathbf{y};M)_j$
are the \textbf{annihilator numbers} of $M$ with respect to the sequence $\mathbf{y}$.
\end{defi}
Clearly, if $\mathbf{y}$ is a filter regular sequence on $M$, then for each $i$ one has that 
$\alpha_{ij}(\mathbf{y};M)$ are equal to zero for almost all $j$.

\begin{prop} \label{2} \emph{[\cite{hh}, Proposition 4.3.5]} 
Let $M$ be a finitely generated 
graded $R$-module and $\mathbf{y}$ a sequence of elements in $R_1$. 
The following conditions are equivalent:
\begin{enumerate}
\item[1.]{$\mathbf{y}$ is a filter regular sequence on $M$;}

\item[2.]{$H_j(y_1,\ldots,y_i;M)$ has finite length for all $j>0$ and all $i$;}

\item[3.]{$H_1(y_1,\ldots,y_i;M)$ has finite length for all $i$.}
\end{enumerate}
Here $H_j(y_1,\ldots,y_i;M)$ denotes the $j$-th homology module of the Koszul complex 
$K.(y_1,\ldots,y_i;M)$.
\end{prop}

In the case $M=R/I$ and the sequence $\mathbf{x} = x_n , x_{n-1},
\ldots,x_1$ one may define some useful annihilator modules. Let $I_0 = I$ 
and $I_i = I_{i-1} + (x_{n-i+1})$ for all $i \in \{1,\ldots,n\}$, one defines
\[
a_{\mathbf{x}}^i(I) = \frac{I_{i-1} : (x_{n-i+1})}{I_{i-1}} .
\]
We remark that $a_{\mathbf{x}}^i(I)$ have finite length for all $i$ if
and only if $\mathbf{x}$ is a filter regular sequence on $R/I$.

\begin{lemma} \label{uso} Let $S$ be a graded $K$-algebra and let $\mathbf{m}$
be the irrelevant ideal. Let $x$ be a homogeneous element of $S$. Then
$(0:x)_j = 0$, for  $j \gg 0$, if and only if $x \notin p$
for all $p \in \emph{Ass}(S), p \ne \mathbf{m}$.
\end{lemma}
\textit{Proof.} We first suppose that $(0:x)$ has finite length and suppose
by contradiction that there exist a relevant associated prime $p$ such that
$x \in p$ and an element $g \in S_i$ such that $(0 : g) = p$. Since $p \neq 
\mathbf{m}$, there exists an element $x \in \mathbf{m} - p$. Now since
$p$ is a prime ideal, we may assume that for all integer $k$, $x^k \notin
(0 : g)$. Further, for any integer $\nu$, we may choose $k$ sufficiently large.
If we set $j = k \cdot \textrm{deg}(x) + i$, then $j > \nu$ and so 
$0 \neq x^kg \in (0:x)_j$, since $x^k \notin (0 : g)$ and $x \in p$. Hence
for any integer $\nu$ there exist elements of $(0:x)$ of higher degree,
a contradiction. Conversely, it is sufficient to show that every element
of $(0:x)$ is nilpotent. In fact we deduce that in high degree (for example
higher then the product of nilpotent orders of a finite system of generators)
there is no element different from zero. Then let $p_i$ for $i = 1,\ldots,k$
the associated primes of $S$ which are not in $\mathbf{m}$ and $q_i$ for $i = 1,\ldots,k$ the related
primary components. Let $J$ the primary component associated to
the maximal ideal. If $y$ is an element of $(0:x)$, then
\[
yx = 0 = \bigcap_{i=1}^{k}q_i \cap J.
\]
So we have $y \in q_i$ for all $i=1,\ldots,k$, otherwise we have that
$x^n \in q_i$, that is $x \in p_i$, a contradiction. We may suppose
$y$ a homogeneous element, in other words $y \in \mathbf{m}$. Then
there exist an integer $r$ such that $y^r \in \mathbf{m}^r \subseteq J$.
Hence $y^r \in \displaystyle\bigcap_{i=1}^kq_i \cap J$, that is $y^r = 0$. $\square$

We use this Lemma to prove the following result:

\begin{prop} \label{ultima} Let $I$ be a monomial ideal in $R$ and let $l(\cdot)$ 
the length function. The following are equivalent
\begin{enumerate}
\item[1.]{$l(a_{\mathbf{x}}^i(I)) < \infty$ for all $i$;}

\item[2.]{$I$ is an ideal of Borel-type.}
\end{enumerate}
\end{prop}
\textit{Proof.} The implication $(2) \Rightarrow (1)$ was proved in 
[\cite{hh}, Proposition 4.3.3]: in fact since $l(a_{\mathbf{x}}^i(I))$ are finite for all $i$,
$\mathbf{x}$ is a filter regular sequence on $R/I$. So using an induction argument and 
Lemma \ref{uso} we are done. Conversely, if $n=2$ there is nothing to prove. Now 
consider the case $n > 2$. We use Proposition \ref{Bti}(2). We may assume that $I+(x_n) / (x_n)$
is of Borel type by induction hypothesis. Hence for any monomial generator $u$ which is not divided
by $x_n$, if $x_i^s | u$, $1 \leq i < j \leq n-1$, there exists $t \geq 0$ such that $x_j^t(u/x_i^s) \in I + (x_n)$
(and hence $x_j^t(u/x_i^s) \in I$). Now suppose $x_n | u = x_1^{s_1} \cdots x_n^{s_n}$, $u \in I$, that is,
$s = s_n \geq 1$. Then we first show that for any $i$ with $1 \leq i \leq n-1$, there exists $k \geq 0$ such that
$x_1^{s_1} \cdots x_i^{s_i+k} \cdots x_{n-1}^{s_{n-1}} \in I$ by induction on $s$. (Indeed, if $s=1$, the
assertion follows from $l((I:x_n)/I) < \infty$). This implies that for any $1 \leq i < j \leq n-1$, there exists
$k_i \geq 0$ such that $x_1^{s_1} \cdots x_i^{s_i+k_i} \cdots \widehat{x_j^{s_j}} \cdots x_{n-1}^{s_{n-1}} \in I + (x_n)$.
In particular $x_1^{s_1} \cdots x_i^{s_i+k_i} \cdots \widehat{x_j^{s_j}} \cdots x_{n-1}^{s_{n-1}} x_{n}^{s_{n}} \in I$, as required.


\section{Preserving extremal Betti numbers and annihilator numbers}
In this section we deal with the annihilator numbers of a graded
$K$-algebra and the correspondence with extremal Betti numbers.
We will use this approach to prove Theorem \ref{principale} based on the Theorem on extremal Betti 
numbers by Bayer, Charalambous, Popescu [see \cite{bcp}].
This section is based on \cite{hh}. The following fundamental lemma will be
useful later.

\begin{lemma} \label{3} Let $I \subseteq R$ be a homogeneous ideal such that
$l(a_{\underline{x}}^i(I))$ are finite for all $i$. Let $<$ be the reverse
lexicographic order. Then:
\[
\dim_K A_{i-1}\left(x_n,x_{n-1},\ldots,x_1 ; \frac{R}{I}\right)_j =
\dim_K A_{i-1}\left(x_n,x_{n-1},\ldots,x_1 ; \frac{R}{\text{in}_<(I)}\right)_j ,
\]
for all $i,j$.
\end{lemma}
\textsc{Proof:} It is suffices to show that the two 
modules above have the same Hilbert function. 
By properties of reverse lexicographic
order, we have that the modules
\[
\frac{R}{(I,x_n,x_{n-1},\ldots,x_{n-i+1})} \ \ , 
\ \ \frac{R}{(\text{in}_<(I),x_n,x_{n-1},\ldots,x_{n-i+1})}
\]
have the same Hilbert function. Consider now the two exact sequences
\[
\begin{CD}
0 \to  A_{i-1}(x_n,\ldots,x_1;R/J)
\to \frac{R}{(J,x_n,x_{n-1},\ldots,x_{n-i+1})}(-1) 
\end{CD}
\]
\[
\begin{CD}
@>\cdot x_{n-i}>>  
\frac{R}{(J,x_n,x_{n-1},\ldots,x_{n-i+1})}  
\to \frac{R}{(J,x_n,x_{n-1},\ldots,x_{n-i})}  \to 0 ,
\end{CD}
\]
for $i=0,\ldots,n$ and $J =I$ or $J =$ in$_<(I)$.
Since the Hilbert function is additive on short exact sequences, we have that
the Hilbert function of $A_{i-1}(x_n,\ldots,x_1;R/I)$ is determined
by the Hilbert function of the modules $R/(I,x_n,x_{n-1},\ldots,x_{n-i})$
and  $R/(I,x_n,x_{n-1},\ldots,x_{n-i+1})$.
This two modules have the same Hilbert function respectively of
$R/(\text{in}_<(I),x_n,x_{n-1},\ldots,x_{n-i})$ and
$R/(\text{in}_<(I),x_n,x_{n-1},\ldots,x_{n-i+1})$, that determine the Hilbert function of the module 
$A_{i-1}(x_n,\ldots,x_1;R/\text{in}_<(I))$. \ $\square$
\ \\
\ \\
Annihilator numbers can also be defined for modules. To do it one just has to 
extend the concept of generic initial ideals to generic initial submodules 
[for details see \cite{e}] and one can show that Theorem \ref{3} holds in the 
general case of a finitely generated graded $R$-module. 

We write $\alpha_{ij}(R/I)$ instead of 
$\alpha_{ij}(x_n, x_{n-1}, \ldots, x_1 ; R/I)$, which are the annihilator
numbers on $R/I$ with respect $x_n, x_{n-1}, \ldots, x_1$.

\begin{coro} \label{4} \upshape
$\alpha_{ij}(R/I) = \alpha_{ij}(R/\text{in}_<(I))$, for all $i,j$. 
\end{coro}

Annihilator numbers of a filtered regular sequence and Betti numbers 
are related to each other. We shall use the convention that
\[
\binom{i}{-1} =
\begin{cases}
\ 0 & \text{if $i \ne -1$} , \\
\ 1 & \text{if $i = -1$} .
\end{cases}
\]

\begin{prop} \label{5} \emph{[\cite{hh}, Proposition 4.3.12]} 
Let $M$ be a finitely generated graded $R$-module
and $\mathbf{y} = y_1,\ldots,y_n$ a $K$-basis of $R_1$ which is a filter
regular sequence on $M$. Then
\[
\beta_{i,i+j}(M) \leqslant \sum_{k=0}^{n-i}\binom{n-k-1}{i-1}
\alpha_{kj}(\mathbf{y};M),
\]
for all $i \geqslant 0$ and all $j$.
\end{prop}

\begin{defi} \upshape Let $M$ be a finitely generated graded $R$-module
and let $\mathbf{y}$ be a $K$-basis of $R_1$ which is a filter regular sequence
on $M$. Let $\alpha_{ij}$ be the annihilator numbers of $M$ with respect
$\mathbf{y}$ and $\beta_{ij}$ be the graded Betti numbers of $M$.
\begin{enumerate}
\item[a.]{An annihilator number $\alpha_{ij} \ne 0$ is called
\textbf{extremal} if $\alpha_{kl} = 0$ for all $(k,l) \ne (i,j)$ with
$k \leqslant i$ and $l \geqslant j$;}

\item[b.]{A graded Betti number $\beta_{i,i+j} \ne 0$ is called
\textbf{extremal} if $\beta_{k,k+l} = 0$ for all $(k,l) \ne (i,j)$ with
$k \geqslant i$ and $l \geqslant j$.}
\end{enumerate}
\end{defi}

Using Proposition \ref{2} and Proposition \ref{5} one can prove 
the following result:
\begin{prop} \label{6} \emph{[\cite{hh}, Theorem 4.3.15]} 
Let $M$ be a graded 
$R$-module and let $\mathbf{y}$ be a $K$-basis of $R_1$ which is a filter 
regular sequence on $M$. Let $\alpha_{ij}$ be the annihilator numbers of $M$ 
with respect to $\mathbf{y}$ and $\beta_{ij}$ be the graded Betti numbers of
$M$. Then $\beta_{i,i+j}$ is an extremal Betti number of $M$ if and only if 
$\alpha_{n-i,j}$ is an extremal annihilator number of $M$. Moreover, if 
the equivalent conditions hold, then
\[
\beta_{i,i+j} = \alpha_{n-i,j}.
\]
\end{prop}
Combining Proposition \ref{6} and Corollary \ref{4} we immediately
obtain the following:
\begin{theo} \label{principale} Let $I \subseteq R$ be a graded ideal. 
Suppose that $l(a_{\mathbf{x}}^i(I))$ are finite for all $i$. Let $<$ be
the reverse lexicographic order.
Then for any two numbers $i,j \in \mathbb{N}$ one has:
\begin{enumerate}
\item[1.]{$\beta_{i,i+j}(I)$ is extremal if and only if $\beta_{i,i+j}(\emph{in}_<(I))$ 
is extremal;}
\item[2.]{if $\beta_{i,i+j}(I)$ is extremal, then $\beta_{i,i+j}(I) = 
\beta_{i,i+j}(\emph{in}_<(I))$.}
\end{enumerate}
\end{theo}

The Theorem above and the results in this section and in section $2$ yield
two results that recover what has been done respectively by Bayer and Stillman in
generic coordinates in \cite{bs} and Trung in \cite{t}:
\begin{coro} Let $I \subseteq R$ be a graded ideal such that
$l(a_{\mathbf{x}}^i(I))$ are finite for all $i$. Let $<$ be the
reverse lexicographic order.
Then
\begin{enumerate}
\item[1.]{\emph{pd}$(I) =$ \emph{pd}$(\emph{in}_<(I))$;}
\item[2.]{\emph{depth}$(R/I) =$ \emph{depth}$(R/\emph{in}_<(I))$;}
\item[3.]{$R/I$ is Cohen-Macaulay if and only if $R/\emph{in}_<(I)$ is 
Cohen-Macaulay;}
\item[4.]{\emph{reg}$(I) =$ \emph{reg}$(\emph{in}_<(I))$.}
\end{enumerate}
\end{coro}

\begin{coro} \emph{[\cite{t}, Theorem 1.3.]} Let $I \subseteq R$ be a graded ideal such that
$l(a_{\mathbf{x}}^i(I))$ are finite for all $i$. Let $<$ be the
reverse lexicographic order.
Then
\begin{enumerate}
\item[1.]{$\emph{reg}_t(I) = \emph{reg}_t(\text{in}_<(I))$;}
\item[2.]{$\emph{a}^{\ast}_t(I) = \emph{a}^{\ast}_t(\text{in}_<(I))$;}
\end{enumerate}
\end{coro}

We conclude this section by investigating on a question about the rigidity of
resolutions in a specific case. Conca in \cite{co1} raised the following
question: let $I$ be a graded ideal and let Gin$(I)$ be his
generic initial ideal with respect the reverse lexicographic order. If
$\beta_i(I) = \beta_i(\text{Gin}(I))$ for some $i$, then it is true that 
$\beta_{k}(I) = \beta_k(\text{Gin}(I))$ for all $k \geqslant i$? This question has
a positive answer as proved by Conca, Herzog and Hibi 
[\cite{chh}, Corollary 2.4]. For $i=0$ this fact was first proved by Aramova,
Herzog and Hibi in \cite{ahh}. 

Consider now $I$ a graded ideal of $R$. We define $I_{\langle j \rangle}$
the ideal generated by all homogeneous polynomial of degree $j$ belonging
to $I$. 
\begin{defi} \upshape A homogeneous ideal $I \subseteq R$ is
\textbf{componentwise linear} (\cite{hh0}), if $I_{\langle j \rangle}$ has a linear 
resolution for all $j$. 
\end{defi}

It is known that $I$ is a componentwise linear ideal in $R$
if and only if $\beta_{ij}(R/I) = \beta_{ij}(R/\text{gin}_<I)$ for all $i,j$ [\cite{ahh}, Theorem 1.1]. 
Here we show an example of an ideal with initial ideal of Borel-type, 
$\mu(I) = \mu(\text{in}_<(I))$ and with different resolution with respect his initial
ideal. We use CoCoA for computations.

\begin{exa} \upshape Let $R = K[x_1,x_2,x_3]$ be the polynomial ring in $3$
variables with the reverse lexicographic order and
\[
\begin{aligned}
I = & ( (2x_1+x_2)^3 , (x_2+2x_3)^3 , (3x_1+x_3)^3  , (x_1+3x_3)^3 , \\
& (3x_1+2x_3)^3 , (2x_2-3x_3)^3 , (4x_1+3x_2)^3 , (3x_1-5x_3)^3 ) .
\end{aligned}
\]
Then $\text{in}_<(I) = (x_1^3 , x_1^2x_2 , x_2^3 , x_1^2x_3 , x_1x_3^2 , 
x_3^3 , x_2^2x_3 , x_1x_2^2).$ It is easy to see that
\[
\text{Ass}\left( \frac{R}{\text{in}_<I} \right) = \{ (x_1 , x_2 , x_3) \}
\]
and so we can conclude that in$_<(I)$ is an ideal of Borel-type.
Further gin$_<(I) = (x_1^3, x_1^2x_2, x_1x_2^2, x_2^3, x_1^2x_3, x_1x_2x_3, 
x_2^2x_3, x_1x_3^2, x_2x_3^3, x_3^4)$.
The Betti tables of the three ideals are the following: \\
\[
\begin{array}{lllll}
\verb"BettiDiagram(I);"\qquad  & \qquad \verb"BettiDiagram(LT(I));" \\
\verb"        0    1    2"   \qquad   & \qquad \verb"        0    1    2"          \\
\verb"--------------------" \qquad  & \qquad  \verb"--------------------"        \\
\verb" 3:     8    9    1"  \qquad   & \qquad \verb" 3:     8    9    2"         \\
\verb" 4:     -    1    2"   \qquad   & \qquad \verb" 4:     -    2    2"          \\
\verb"--------------------"  \qquad  & \qquad  \verb"--------------------"        \\
\verb"Tot:    8   10    3" \qquad & \qquad \verb"Tot:    8   11    4" 
\end{array}
\]
\begin{center}
\verb"BettiDiagram(Gin(I));"\\
\verb"        0    1    2"\\
\verb"--------------------"\\
\verb" 3:     8   11    4"\\
\verb" 4:     2    4    2"\\
\verb"--------------------"\\
\verb"Tot:   10   15    6"\\
\end{center}
We immediately notice that $\mu(I) = \mu(\text{in}_<(I)) \neq \mu(\text{gin}_<(I))$,
$\beta_{i}(I) \neq \beta_i(\text{in}_<(I))$ for all $i = 1,2$, and
$\beta_i(I) \neq \beta_i(\text{gin}_<(I))$ for all $i$. Then this is an example of
a monomial ideal, in$_<(I)$, that is Borel-type but not componentwise linear. In
particular if $I$ is a graded ideal with initial ideal of Borel-type such that $\mu(I) =
\mu(\text{in}_<(I))$, it is not true in general that all the $\beta_i$'s are equal.
\end{exa}

\begin{oss} \upshape Conca in a private communication reported that if $I$ is a graded 
ideal such that in$_<(I)$ is componentwise linear and $\mu(I) = \mu(\text{in}_<(I))$, then 
$\beta_i(I) = \beta_i(\text{in}_<(I))$ for all $i$.
\end{oss}

\section{Preserving reduction number}
Let $R = K[x_1,\ldots,x_n]$ be a polynomial ring, $K$ an infinite field and $I \ne 0$
a homogeneous ideal in $R$. We set $\mathbf{m} = (x_1,\ldots,x_n)R/I$.
A homogeneous ideal $J \subset \mathbf{m}$ is called
a reduction of $\mathbf{m}$ if $\mathbf{m}^{r+1} = J \mathbf{m}^r$ for some integer
$r \geqslant 0$. $J$ is called a minimal reduction, if it is minimal with respect to 
inclution. The \textbf{reduction number} of $\mathbf{m}$ with
respect to a minimal reduction $J$ of $\mathbf{m}$, denoted by $r_{J}(\mathbf{m})$
or $r_{J}(R/I)$, is the smallest $r \geqslant 0$ such that $\mathbf{m}^{r+1} = 
J \mathbf{m}^r$. The reduction number of $\mathbf{m}$, denoted by
$r(\mathbf{m})$ or $r(R/I)$, is the infimum of $r_{J}(\mathbf{m})$ over all possible
minimal reductions $J$ of $\mathbf{m}$. For the reductions theory see \cite{nr}.
Consider now in$(I)$, the initial ideal of $I$ with respect to some admissible
term order on the terms of $R$. Vasconcelos in \cite{v} conjectured that
\[
r(R / I) \leqslant r(R / \text{in}(I)) .
\]
Bresinsky and Hoa proved, in \cite{bh}, that the conjecture is true for generic coordinates. 
Trung in \cite{t} proved that the equality holds in generic coordinates with 
respect the reverse lexicographic order.
Moreover the conjecture was proved by Conca in \cite{co} and independently by Trung.
We see now that if $I$ is an ideal such that the lengths $l(a_{\mathbf{x}}^i)$ are finite,
the equality of reduction numbers of $I$ and in$(I)$ is not necessarily reached.
So assume that dim$(R/I) = d$ and denote by $\overline{y}$ 
the image of $y \in R$ in $R/I$.

\begin{defi} \upshape If $M$ is any graded $R$-module of finite lenght, we define
\[
a(M) =
\begin{cases}
\max\{p : M_p \ne 0\} & \text{if $M \ne 0$} ; \\
\ \ \ \ \ \ \ \ - \infty & \text{if $M = 0$} .
\end{cases}
\]
\end{defi}

By $\cite{nr}$ a minimal reduction of $\mathbf{m}$ can always be generated by
a system of parameters. Furthermore worth the following

\begin{lemma} \label{lemma3} \emph{[\cite{bh}, Lemma 3]} The ideal $J = (\overline{y}_1 , \ldots , 
\overline{y}_d) \subseteq R/I$ is a minimal reduction of $\mathbf{m}$
if and only if $\{ \overline{y}_1 , \ldots , \overline{y}_d \}$ is a system of
parameters (s.o.p.) of $R/I$ with $\overline{y}_i$ linear forms, 
$1 \leqslant i \leqslant d$ . In this case
\[
r_J(R/I) = a(R/(I,y_1,\ldots,y_d)) .
\]
\end{lemma}

\begin{prop} \label{prop4} \emph{[\cite{bh}, Proposition 4]} Let $I$ be a homogeneous ideal
in $R = K[x_1,\ldots,x_n]$. Then,
\[ 
r(R/I) \geq \emph{min}\{ \emph{deg}(F) : F \in I, \ F \ homogeneous \} - 1 .
\]
\end{prop}

\begin{lemma} \label{lemma5} \emph{[\cite{bh}, Lemma 5]} Assume that $I$ is a monomial
ideal of $R$ such that $\overline{x}_{n-d+1} , \ldots , \overline{x}_n$ is a s.o.p. of $R/I$. 
Then any minimal reduction $J$ of $\mathbf{m}$ is generated by $d$ linear forms 
$\overline{y}_1, \ldots , \overline{y}_d$ with
\[
y_i = x_{n-d+1} + a_{i,1}x_1 + \cdots + a_{i,n-d}x_{n-d}, \ 1 \leq i \leq d,
\] 
where $a_{i,j} \in K$, $1 \leq j \leq n-d$.
\end{lemma}

\begin{coro} \label{coro6} \emph{[\cite{bh}, Corollary 6]} Assume that $I$ is a monomial
ideal of $R$ such that $\overline{x}_{n-d+1} , \ldots , \overline{x}_n$ is a s.o.p. of $R/I$.
Then for any minimal reduction $J$ of $\mathbf{m}$ we have
\[
r_J(R/I) \leq r_{(x_{n-d+1},\ldots,x_n)}(R/I) .
\]
\end{coro}

We notice that an ideal $I$ such that the lengths $l(a^i_{\mathbf{x}})$ are finite for $i = 1 , \ldots , d$,
satisfies the hypothesis of Corollary \ref{coro6} but in general the equality may not be
reached:
\begin{exa} \upshape \label{exx} [\cite{bh}, Example 7] Consider the ideal $I = (x_1^4, x_1x_2^3,x_1x_3^2)$ in
$R = K[x_1,x_2,x_3]$. It is easy to see that $I$ is of Borel-type and $d =$ dim$(R/I) = 2$. Equivalently
$\{x_2,x_3\}$ is a filter regular sequence in $R/I$ and so a s.o.p in $R/I$. Hence $(x_2,x_3)$ is
a minimal reduction of $\mathbf{m}$ in $R/I$ by Lemma \ref{lemma3}. Again by Lemma \ref{lemma3} we get
\[r_J(R/I) =
\begin{cases}
\ \ 3 & \text{if} \ J = (x_2,x_3), \\
\ \ 2 & \text{if} \ J = (x_2,x_3-x_1) .
\end{cases}
\]
Since $3$ is the least degree in the generating set of $I$, by Proposition \ref{prop4},
we have that $r(R/I) = 2 < r_{(x_2,x_3)}(R/I)$.
\end{exa}

The following lemma generalize Lemma 4.1 in \cite{t}.

\begin{lemma} \label{mylemma} Let \emph{in}$(I)$ be the initial ideal of $I$ with respect revlex. 
Let $d =$ \emph{dim}$(R/I)$. Let $J = (I,x_{n-d+1},\ldots,x_n)/I$ and $K = (\emph{in}(I),x_{n-d+1},\ldots,x_n)/\emph{in}(I)$. If $l(a_{\mathbf{x}}^i)$ are finite for all $i \in \{1,\ldots,d\}$, then $r_J(R/I) = r_K(R/\emph{in}(I))$.
\end{lemma}
\textsc{Proof:} By hypothesis $\overline{x}_{n-d+1} , \ldots , \overline{x}_n$ is a filter regular sequence in $R/I$
and so it is a s.o.p in $R/I$. Then, using Lemma \ref{lemma3}, $J$ is a minimal reduction of 
$\mathbf{m}$ in $R/I$ such that $r_J(R/I) = a(R/(I,x_{n-d+1},\ldots,x_n))$. Since $R/I$ 
and $R/\text{in}(I)$ share the same Hilbert function, $\overline{x}_{n-d+1} , \ldots , \overline{x}_n$ is
also a s.o.p in $R/\text{in}(I)$ and so by Lemma \ref{lemma3},
$K$ is a minimal reduction of the homogeneous maximal ideal in $R/\text{in}(I)$ such that
$r_K(R/\text{in}(I)) = a(R/(\text{in}(I),x_{n-d+1},\ldots,x_n))$. Since we use the reverse lexicographic order,
we have the following identity:
\[
\text{in}(I,x_{n-d+1},\ldots,x_n) = (\text{in}(I),x_{n-d+1},\ldots,x_n).
\]
Hence $a(R/(\text{in}(I),x_{n-d+1},\ldots,x_n)) = a(R/\text{in}(I,x_{n-d+1},\ldots,x_n))$. Finally,
since $R/I$ and $R/\text{in}(I)$ share the same Hilbert function, we obtain that
$a(R/(\text{in}(I),x_{n-d+1},\ldots,x_n)) = a(R/(I,x_{n-d+1},\ldots,x_n))$, as required. $\square$

\begin{oss} \upshape Let $I$ be an ideal such that the lengths $l(a^i_{\mathbf{x}})$ 
are finite for all $i = 1 , \ldots , d$ and let in$(I)$ be the initial ideal of $I$ with respect the
reverse lexicographic order. Under these hypothesis, we observe that the reduction numbers of
$I$ and in$(I)$ are not necessarily equal. In fact, by hypothesis, $x_n , \ldots , x_{n-d+1}$ is a filter 
regular sequence in $R/I$ and so a s.o.p in $R/I$. Hence the ideal $K = (I, x_{n-d+1} , \ldots , x_n)/I$ 
is a minimal reduction of $\mathbf{m}$ in $R/I$. We first suppose that for all minimal reductions 
$J$ in $R/I$, $r(R/I) = r_K(R/I) \leq r_J(R/I)$. Using Lemma \ref{mylemma} we obtain that
\[
r(R/I) = r_K(R/I) = r_{K'}(R/\text{in}(I)) \geq r(R/I),
\]
where $K' = (\text{in}(I),x_{n-d+1},\ldots,x_n)/\text{in}(I)$ is a minimal reduction in $R/\text{in}(I)$.
Since $r(R/I) \leq r(R/\text{in}(I))$ is true in general, we obtain that $r(R/I) = r(R/\text{in}(I))$.
By Example \ref{exx}, 
it might exist a minimal reduction $J$ in $R/I$ such that
$r(R/I) = r_J(R/I) < r_K(R/I)$. By Corollary \ref{coro6}, we know that
$r_{J'}(R/\text{in}(I)) \leq r_{K'}(R/\text{in}(I))$ for all minimal reductions $J'$ in $R/\text{in}(I)$.
In particular, we can pick $J'$ such that $r(R/\text{in}(I)) = r_{J'}(R/\text{in}(I))$.
Suppose now that $r_{J'}(R/\text{in}(I)) = r_{K'}(R/\text{in}(I))$. 
In this case using Lemma \ref{mylemma} we have
\[
r(R/I) < r_K(R/I) = r_{K'}(R/\text{in}(I)) = r(R/\text{in}(I)),
\]
that is $r(R/I) < r(R/\text{in}(I))$. Conversely suppose
$r_{J'}(R/\text{in}(I)) < r_{K'}(R/\text{in}(I))$. In this case nothing
can be said more than just the well known inequality $r(R/I) \leq r(R/\text{in}(I))$.
\end{oss}

\section{Quasi-stable versus Borel type ideals and Pommaret bases}

Let $R = K[x_1,\ldots,x_n]$ be the polynomial
ring over a field $K$ in $n$ variables and 
$\mu = [\mu_1 , \ldots , \mu_n]$ be an exponent vector, with $x^{\mu}$
we denote a monomial in $R$ and with $f$ a polynomial such that
in$_<(f) = x^{\mu}$ with respect to the reverse lexicographic order. 
The following definitions and Proposition \ref{QS} are in \cite{hss}
and \cite{s}.

\begin{defi} \upshape We define the \textbf{class} of $\mu$ as the
integer
\[
\text{cls}(\mu) = \text{min}\{i : \mu_i \neq 0\}.
\]
If $f$ is a polynomial with in$_<(f) = x^{\mu}$, by cls$(f)$ one means cls$(\mu)$.
Then the \textbf{multiplicative variables} of $f \in R$ (or $x^{\mu}$) are
\[
\mathrm{X}_R(f) = \mathrm{X}_R(x^{\mu}) = \{x_1,\ldots,x_{\text{cls}\mu}\}.
\]
\end{defi}
If we consider $f = x_2^2x_3 + x_3x_4^2+x_3^3 \in R = K[x_1,x_2,x_3,x_4]$, then we
have that $\mathrm{X}_R(f) = \{x_1,x_2\}$, since \text{cls}$(\text{in}_<(f)) =$
\text{cls}$(x_2^2x_3) = 2$.

\begin{defi} \upshape We say that $x^{\mu}$ is an \textbf{involutive divisor}
of $x^{\nu}$, with $\nu$ another index vector, if $x^{\mu} | x^{\nu}$ and
$x^{\nu-\mu} \in K[\mathrm{X_R(x^{\mu})}]$.
\end{defi}
For example $x_2x_3^2$ is not an involutive divisor of $x_2^2x_3^3$, since
his class is two and $x_2x_3 \notin K[x_1,x_2]$. Instead $x_2x_3^3$ is
an involutive divisor of $x_2^2x_3^3$, since his class is two and $x_2 \in 
K[x_1,x_2]$. In the same way one can see that also $x_3^3$ is an involutive
divisor of $x_2^2x_3^3$.

\begin{defi} \upshape Let $\mathcal{H} \subseteq R$ be a finite subset of only terms.
We say that $\mathcal{H}$ is a \textbf{Pommaret basis} of the monomial ideal
$I = (\mathcal{H})$, if
\[
\bigoplus_{h \in \mathcal{H}}K[\mathrm{X_R(h)}] \cdot h = I
\]
as $K$-linear space. (In this case each term $x^{\nu} \in I$ has a unique involutive
divisor $x^{\mu} \in \mathcal{H}$)
A finite polynomial set $\mathcal{H} \subseteq R$ is a Pommaret basis of a polynomial
ideal $I$ with respect a monomial order $<$ (revlex in our case), if all elements of 
$\mathcal{H}$ have distinct leading terms and these terms form a Pommaret
basis of the ideal \text{lt}$_<(I)$.
\end{defi}

\begin{defi} \upshape A monomial ideal $I$ is \textbf{quasi-stable},
if possess a finite Pommaret basis.
\end{defi}

Similarly to the Borel type case, we give here a characterization 
of quasi-stable ideal:

\begin{prop} \label{QS} \emph{[\cite{hss}, Theorem 11]} Let $I \subseteq R$ be a
monomial ideal and $d =$ \emph{dim}$(R/I)$. Then the following
are equivalent.
\begin{enumerate}
\item[1.]{$I$ is quasi-stable;}

\item[2.]{$x_1$ is a non zero divisor for $R/I^{\text{sat}}$, where $I^{\text{sat}} = 
I : \mathbf{m}^k$ is the \textbf{saturation} of $I$. Besides, for all $1 \leqslant k 
< d$, $x_{k+1}$ is a non zero divisor for $R/(I,x_1,\ldots,x_k)^{\text{sat}}$;}

\item[3.]{$I : x_1^\infty \subseteq I : x_2^\infty \subseteq \cdots \subseteq 
I : x_d^\infty$ and for all $d < k \leqslant n$, there exists an exponent 
$e_k \geqslant 1$ such that $x_k^{e_k} \in I$;}

\item[4.]{for all $1 \leqslant k \leqslant n$ the equality $I : x_k^\infty = 
I : (x_k,\ldots,x_n)^\infty$ holds;}

\item[5.]{if $p \in$ Ass$(R/I)$, then $p = (x_k,\ldots,x_n)$ for some $k$;}

\item[6.]{if $x^{\mu} \in I$ and $\mu_i > 0$ for some $1 \leqslant i < n$, then
for some $0 < r \leqslant \mu_i$ and $i < j \leqslant n$, there exists an integer $s \geqslant 0$
such that $x_j^sx^{\mu}/x_i^r \in I$.}
\end{enumerate}
\end{prop}
\textit{Proof:} for the proof see [Proposition 4.4, \cite{s}]. $\square$

\begin{exa} \upshape Consider the ideal $I = (x_1x_3 , x_2x_3 , x_3^2)$ in
$K[x_1,x_2,x_3]$. We claim that $\mathcal{H} = [x_1x_3 , x_2x_3 , x_3^2]$
is a Pommaret basis of $I$. In fact one can easily see that $I = (x_1,x_2,x_3^2) \cap (x_3)$
satisfies the condition $5$ of Proposition \ref{QS} (see also Example \ref{exfin}).

However $I$ is not an ideal of Borel type. In fact 
using $(3)$ in Proposition \ref{Bti} we have that
$x_3 | x_1x_3$ but there exists no integer $t \geqslant 0$ such that 
$x_1^{t+1} \in I$, otherwise we have
\[
x_1^{t+1} = x_1x_3 f + x_2x_3 g + x_3^2 h
\]
for some $f,g,h \in R$ and so $x_3 | x_1^{t+1}$, a
contradiction.
\end{exa}

\begin{exa} \label{exfin} \upshape Let $R = K[x_1,x_2,x_3]$.
If we consider $I = (x_1x_3 , x_2x_3 , x_3^2)$ the ideal of
the previous example, we have that
\[
\text{Ass}(R/I) = \left \{ (x_3) , (x_1,x_2,x_3) \right \}.
\]
All the associated ideal of $I$ are of the form $(x_k,\ldots,x_3)$ for some 
$k \leqslant 3$ ($k=3$ and $k=1$). Using $(5)$ in Proposition \ref{QS} we have that
$I$ is quasi-stabile but not Borel type, since there exists no $k$ such that 
$(x_3) = (x_1,\ldots,x_k)$. We show now an example of ideal $J$ that is Borel type
but not quasi-stable. We can transform $I$ using the change of variables $x_1 \to x_3$,
$x_2 \to x_2$ and $x_3 \to x_1$ and we obtain the ideal $J = (x_1x_2,x_1x_3,
x_1^2)$. In this case we have
\[
\text{Ass}(R/J) = \left\{ (x_1) , (x_1,x_2,x_3) \right\}
\]
and so the associated primes of $J$ are of the form $(x_1,\ldots,x_k)$ for some
$k \leqslant 3$ ($k=1$ and $k=3$). So $I$ is Borel type but it is not quasi-stabile 
since there exists no $k \leqslant 3$ such that $(x_1) = (x_k,\ldots,x_3)$.
\end{exa}

\section*{Acknowledgment}
We thank Professor Maria Evelina Rossi for her support and guidance
in investigating the subject and the Department of Mathematics of
University of Genova for hospitality.

\newpage


\end{document}